\documentclass{article}

\usepackage{amsthm}
\usepackage{amssymb}
\usepackage{amsxtra}

\newtheorem{thm}{Theorem}[section]

\newtheorem{conj}{Conjecture}[section]
\newtheorem{question}{Question}

\title{Finding Almost Squares II}
\author{Tsz Ho Chan}

\begin{document}
\maketitle
\begin{abstract}
In this article, we study short intervals that contain another type of ``almost square", an integer $n$ which can be factored in two different ways $n = a_1 b_1 = a_2 b_2$ with $a_1, a_2, b_1, b_2$ close to $\sqrt{n}$.
\end{abstract}
\section{Introduction}
In [\ref{C}], the author studied the problem of finding ``almost squares" in short intervals, namely:
\begin{question}
\label{old}
For $0 \leq \theta < 1/2$, what is the least $f(\theta)$ such
that, for some $c_1, c_2 > 0$, any interval $[x - c_1 x^{f(\theta)}, x
+ c_1 x^{f(\theta)}]$ contains an integer $n$ with $n = ab$, where
$a$, $b$ are integers in the interval $[x^{1/2} - c_2 x^\theta,
x^{1/2} + c_2 x^\theta]$? Note: $c_1$ and $c_2$ may depend on
$\theta$.
\end{question}
A similar question is the following:
\begin{question}
\label{new}
For $0 \leq \theta < 1/2$, what is the least $g(\theta)$ such
that, for some $c_1, c_2 > 0$, any interval $[x - c_1 x^{g(\theta)}, x
+ c_1 x^{g(\theta)}]$ contains an integer $n$ with $n = a_1 b_1 = a_2 b_2$, where
$a_1 < a_2 \leq b_2 < b_1$ are integers in the interval $[x^{1/2} - c_2 x^\theta, x^{1/2} + c_2 x^\theta]$? Note: $c_1$ and $c_2$ may depend on
$\theta$.
\end{question}
Note: Actually, the author first considered Question \ref{new} and then turned to Question \ref{old} which has connection to the problem on the distribution of $n^2 \alpha \pmod 1$ and the problem on gaps between sums of two squares.

In [\ref{C}], we showed that $f(\theta) = 1/2$ when $0 \leq \theta < 1/4$, $f(1/4) = 1/4$ and $f(\theta) \geq 1/2 - \theta$. We conjectured that $f(\theta) = 1/2 - \theta$ for $1/4 < \theta < 1/2$ and gave conditional result when $1/4 < \theta < 3/10$. To Question \ref{new}, we have the following
\begin{thm}
\label{theorem1}
For $0 < \theta < 1/4$, $g(\theta)$ does not exist (i.e. all possible products of pairs of integers in $[x^{1/2} - c_2 x^\theta, x^{1/2} + c_2 x^\theta]$ are necessarily distinct for large $x$).
\end{thm}
\begin{thm}
\label{theorem2}
For $1/4 \leq \theta < 1/2$, $g(\theta) \geq 1 - 2\theta$.
\end{thm}
\begin{thm}
\label{theorem3}
For $1/4 \leq \theta \leq 1/3$, $g(\theta) \leq 1 - \theta$.
\end{thm}
We believe that the lower bound is closer to the truth and conjecture
\begin{conj}
\label{conj1}
For $1/4 \leq \theta < 1/2$, $g(\theta) = 1 - 2\theta$.
\end{conj}
\section{Preliminaries and $0 \leq \theta < 1/4$}
Suppose $n = a_1 b_1 = a_2 b_2$ with $x^{1/2} - c_2 x^\theta \leq a_1 < a_2 \leq b_2 < b_1 \leq x^{1/2} + c_2 x^\theta$. Let $d_1 = (a_1, a_2)$ and $d_2 = (b_1, b_2)$ be the greatest common divisors. Then we must have $d_1, d_2 > 1$. For otherwise, say $d_1 = 1$, then $a_2$ divides $b_1$ which implies $x^{1/2} + c_2 x^\theta \geq b_1 \geq 2 a_2 \geq 2 x^{1/2} - 2 c_2 x^\theta$. This is impossible for large $x$ as $\theta < 1/2$. Now, let $a_1 = d_1 e_1$, $a_2 = d_1 e_2$, $b_1 = d_2 f_1$ and $b_2 = d_2 f_2$. Here $(e_1, e_2) = 1 = (f_1, f_2)$. Then
$$n = d_1 e_1 d_2 f_1 = d_1 e_2 d_2 f_2 \mbox{ gives } e_1 f_1 = e_2 f_2.$$
Due to co-primality, $e_2 = f_1$ and $e_1 = f_2$. Therefore,
\begin{equation}
\label{form}
n = (d_1 e_1) (d_2 e_2) = (d_1 e_2) (d_2 e_1)
\end{equation}
with $1< d_1 < d_2$, $e_1 < e_2$ and $(e_1, e_2) = 1$.

\smallskip

Now, from $a_2 - a_1 \leq 2c_2 x^\theta$, $d_1 \leq d_1 e_2 - d_1 e_1 \leq 2c_2 x^\theta$. Similarly, one can deduce that $d_2, e_1, e_2 \leq 2c_2 x^\theta$. Moreover, as $d_1 e_1 = a_1 \geq x^{1/2} - c_2 x^\theta$, we have $d_1, e_1 \geq \frac{1}{2c_2} x^{1/2 - \theta} - \frac{1}{2}$. Similarly, $d_2, e_2 \geq \frac{1}{2c_2} x^{1/2 - \theta} - \frac{1}{2}$. Summing up, we have
\begin{equation}
\label{range}
\frac{1}{2c_2} x^{1/2 - \theta} - \frac{1}{2} \leq d_1, d_2, e_1, e_2 \leq 2c_2 x^\theta.
\end{equation}
From (\ref{range}), we see that no such $n$ exists for $0 \leq \theta < 1/4$ and hence Theorem \ref{theorem1}.
\section{Lower bound for $g(\theta)$}
From (\ref{form}) and (\ref{range}), we see that an integer $n = a_1 b_1 = a_2 b_2$, satisfying the conditions for $a_1, a_2, b_1, b_2$ in Question \ref{new}, must be of the form:
$$n = (d_1 e_1) (d_2 e_2) \mbox{ with } \frac{1}{2c_2} x^{1/2 - \theta} - \frac{1}{2} \leq d_1, d_2, e_1, e_2 \leq 2c_2 x^\theta$$
and $x^{1/2} - c_2 x^\theta \leq d_1 e_1 < d_1 e_2, \, d_2 e_1 < d_2 e_2 \leq x^{1/2} + c_2 x^\theta$. In particular, $e_2 d_2 - e_2 d_1 \leq 2c_2 x^\theta$ which implies $e_2 - e_1 \leq 2c_2 x^\theta / d_2$. Similarly, $d_2 - d_1 \leq 2c_2 x^\theta / e_2$. Thus, the number of such quartuple $(d_1, d_2, e_1, e_2)$ is bounded by
$$\ll \mathop{\sum_{x^{1/2 - \theta} \ll d_2, e_2 \ll x^\theta}}_{x^{1/2} - c_2 x^\theta \leq d_2 e_2 \leq x^{1/2} + c_2 x^\theta} \frac{x^\theta}{e_2} \frac{x^\theta}{d_2} \ll \frac{x^{2\theta}}{x^{1/2}} x^\theta x^\epsilon = x^{3\theta - 1/2 + \epsilon}$$
for any $\epsilon > 0$ as the number of divisor function $d(n) \ll n^\epsilon$. It follows that there are at most $O(x^{3\theta - 1/2 + \epsilon})$ such integers $n$ in the interval $[x - c_2 x^{1/2+\theta}/3, x + c_2 x^{1/2+\theta}/3]$. Therefore, some two consecutive such $n$'s have gap
$$\gg \frac{x^{1/2 + \theta}}{x^{3\theta - 1/2 + \epsilon}} = x^{1 - 2\theta - \epsilon}.$$
Pick $y$ to be the midpoint between these two integers. Then, for some constant $c > 0$, the interval $[y - c y^{1-2\theta-\epsilon}, y + c y^{1-2\theta-\epsilon}]$ does not contain any integer $n = a_1 b_1 = a_2 b_2$ with $y^{1/2} - c_2 y^\theta/2 \leq a_1 < a_2 \leq b_2 < b_1 \leq y^{1/2} + c_2 y^\theta/2$ as $x - c_2 x^{1/2 + \theta}/3 \leq y \leq x + c_2 x^{1/2 + \theta}/3$. Consequently, for any constant $c, c' > 0$, there is arbitrarily large $y$ such that the interval $[y - c y^{1-2\theta-2\epsilon}, y + c y^{1-2\theta-2\epsilon}]$ does not contain any integer $n = a_1 b_1 = a_2 b_2$ with $y^{1/2} - c' y^\theta \leq a_1 < a_2 \leq b_2 < b_1 \leq y^{1/2} + c' y^\theta$. Therefore, $g(\theta) \geq 1 - 2\theta - 2\epsilon$ which gives Theorem \ref{theorem2} by letting $\epsilon \rightarrow 0$.

\section{Upper bound for $g(\theta)$}
Proof of Theorem \ref{theorem3}: For any large $x$, set $N = [x^{1/4}]$ and $\xi = \{x^{1/4}\}$, the integer part and fractional part of $x^{1/4}$ respectively. Based on (\ref{form}), we are going to pick, for $0 \leq \epsilon \leq 1/2$,
\begin{equation}
\label{de}
d_1 = q N + r_1, \; d_2 = q N + r_2, \; e_1 = \frac{N + s_1}{q}, \; e_2 = \frac{N + s_2}{q}
\end{equation}
for some $1 \leq q \leq N^\epsilon$, $0 \leq r_1, r_2 < N$ and $s_1, s_2 \ll q$ with $N \equiv -s_1 \equiv -s_2 \pmod q$. Our goal is to make
$$x = (N+\xi)^4 = N^4 + 4N^3 \xi + O(N^2) \approx (qN+r_1) \frac{N + s_1}{q} (qN+r_2) \frac{N + s_2}{q}.$$
The right hand side above is
\begin{equation}
\label{approx}
\begin{split}
=& \Bigl[N^2 + \Bigl(\frac{r_1}{q} + s_1\Bigr)N + \frac{r_1 s_1}{q} \Bigr] \Bigl[N^2 + \Bigl(\frac{r_2}{q} + s_2\Bigr)N + \frac{r_2 s_2}{q} \Bigr] \\
=& N^4 + \Bigl(\frac{r_1 + r_2}{q} + s_1 + s_2 \Bigr)N^3 + \Bigl[\frac{r_1 s_1}{q} + \frac{r_2 s_2}{q} + \Bigl(\frac{r_1}{q} + s_1\Bigr) \Bigl(\frac{r_2}{q} + s_2\Bigr) \Bigr] N^2 \\
&+ \Bigl[\frac{r_1 s_1}{q} \Bigl(\frac{r_2}{q} + s_2\Bigr) + \frac{r_2 s_2}{q} \Bigl(\frac{r_1}{q} + s_1\Bigr)\Bigr] N + \frac{r_1 s_1 r_2 s_2}{q^2}
\end{split}
\end{equation}
By Dirichlet's Theorem on diophantine approximation, we can find integer $1 \leq q \leq N^\epsilon$ such that
$$\Big| 4\xi - \frac{p}{q} \Big| \leq \frac{1}{q N^\epsilon}$$
for some integer $p$. Fix such a $q$. Then, pick $s_1 < s_2 < 0$ to be the largest two integers such that $N \equiv -s_1 \equiv -s_2 \pmod q$. Clearly, $s_1, s_2 \ll q$. Then, one simply picks some $0 < r_1 < r_2 \ll q^2$ such that $\frac{r_1 + r_2}{q} + s_1 + s_2 = \frac{p}{q}$. With these values for $q, r_1, r_2, s_1, s_2$, (\ref{approx}) is
$$=N^4 + 4N^3 \xi + O(N^{3-\epsilon}) + O(q^2 N^2) + O(q^3 N) + O(q^4).$$
Hence, we have just constructed an integer $n = d_1 e_1 d_2 e_2$ which is within $O(N^{3-\epsilon}) + O(N^{2+2\epsilon}) = O(x^{3/4 - \epsilon/4}) + O(x^{1/2 + \epsilon/2}) = O(x^{3/4 - \epsilon/4})$ from $x$ if $\epsilon \leq 1/3$. One can easily check that $a_1 = d_1 e_1$, $b_1 = d_2 e_2$, $a_2 = d_1 e_2$ and $b_2 = d_2 e_1$ are in the interval $[x^{1/2} - C x^{1/4 + \epsilon/4}, x^{1/2} + C x^{1/4 + \epsilon/4}]$ for some constant $C > 0$. Set $\theta = 1/4 + \epsilon/4$. We have, for some $C' > 0$, $n = a_1 b_1 = a_2 b_2$ in the interval $[x - C'x^{1 - \theta}, x+C'x^{1 - \theta}]$ such that $a_1 < a_2, b_2 < b_1$ are integers in $[x^{1/2} - Cx^\theta, x^{1/2} + Cx^\theta]$ provided $1/4 \leq \theta \leq 1/4 + 1/12 = 1/3$. This proves Theorem \ref{theorem3}.
\section{Open questions}
Conjecture \ref{conj1} may be too hard to prove in the moment. As a starting point, can one show that $g(1/4) = 1/2$? Or even $g(1/4) < 3/4$? Another possibility would be trying to get some results conditionally like [\ref{C}]. Also, one may consider $g(\theta)$ when $\theta$ is near to $1/2$. This leads to the problem about gaps between integers that have more than one representation as a sum of two squares.

\bigskip

{\bf Acknowledgement} The author would like to thank the American Institute of Mathematics for providing a stimulating environment to work at.

Tsz Ho Chan\\
American Institute of Mathematics\\
360 Portage Avenue\\
Palo Alto, CA 94306\\
U.S.A.\\
thchan@aimath.org

\end{document}